\documentclass[a4paper,11pt]{amsart}
\usepackage{amssymb}
\usepackage[arrow,matrix]{xy}

\DeclareMathOperator{\Ext}{Ext}
\DeclareMathOperator{\Hom}{Hom}

\DeclareMathOperator{\GL}{GL}
\DeclareMathOperator{\spec}{spec}

\DeclareMathOperator{\Mod}{mod}
\DeclareMathOperator{\gldim}{gl.\!dim}

\newcommand{\bp}{\mathbf{p}}
\newcommand{\bl}{\pmb{\lambda}}

\newcommand{\gL}{\Lambda}
\newcommand{\cA}{\mathcal{A}}
\newcommand{\cB}{\mathcal{B}}
\newcommand{\cD}{\mathcal{D}}
\newcommand{\cH}{\mathcal{H}}

\newcommand{\bZ}{\mathbb{Z}}
\newcommand{\bR}{\mathbb{R}}
\newcommand{\bform}[2]{\left\langle {#1}, {#2} \right\rangle}
\newcommand{\Bform}{\bform{\cdot}{\cdot}}

\newcommand{\dH}{\cD^b(\cH)}
\newcommand{\dL}{\cD^b(\gL)}
\newcommand{\dX}{\cD^b(X)}

\newcommand{\one}{\mathbf{1}}

\theoremstyle{plain}
\newtheorem{thm}{Theorem}[section]
\newtheorem{lem}[thm]{Lemma}
\newtheorem{prop}[thm]{Proposition}
\newtheorem{cor}[thm]{Corollary}

\theoremstyle{definition}
\newtheorem{defn}[thm]{Definition}
\newtheorem{exmp}[thm]{Example}
\newtheorem{rem}[thm]{Remark}

\numberwithin{equation}{section}

\title[On the periodicity of Coxeter transformations]
{On the periodicity of Coxeter transformations and the non-negativity of their Euler forms}

\author{Sefi Ladkani}
\address{Einstein Institute of Mathematics, The Hebrew University of Jerusalem, Jerusalem 91904, Israel}
\email{sefil@math.huji.ac.il}

\begin{document}

\begin{abstract}
We show that for piecewise hereditary algebras, the periodicity of the
Coxeter transformation implies the non-negativity of the Euler form.
Contrary to previous assumptions, the condition of piecewise heredity
cannot be omitted, even for triangular algebras, as demonstrated by
incidence algebras of posets.

We also give a simple, direct proof, that certain products of
reflections, defined for any square matrix $A$ with 2 on its main
diagonal, and in particular the Coxeter transformation corresponding to
a generalized Cartan matrix, can be expressed as $-A_{+}^{-1} A_{-}^t$,
where $A_{+}$, $A_{-}$ are closely associated with the upper and lower
triangular parts of $A$.
\end{abstract}



\maketitle

\section{Introduction}
Let $V$ be a free abelian group of finite rank and let $\Bform : V
\times V \to \bZ$ be a non-degenerate $\bZ$-bilinear form on $V$.
The \emph{Coxeter transformation} $\Phi : V \to V$ corresponding to
$\Bform$ is defined via the equation $\bform{x}{y} = -
\bform{y}{\Phi x}$ for $x, y \in V$~\cite{Lenzing96}.

The purpose of this paper is to study the relations between positivity
properties of the form $\Bform$ and periodicity properties of its
Coxeter transformation $\Phi$. Recall that $\Bform$ is \emph{positive}
if $\bform{x}{x} > 0$ for all $0 \neq x \in V$, \emph{non-negative} if
$\bform{x}{x} \geq 0$ for all $x \in V$ and \emph{indefinite}
otherwise. The transformation $\Phi$ is \emph{periodic} if $\Phi^m$
equals the identity $I$ for some integer $m \geq 1$ and \emph{weakly
periodic} \cite{Sato05} if $(\Phi^m-I)^n=0$ for some integers $m,n \geq
1$.

Implications in one direction are given in the paper~\cite{Sato05},
where linear algebra techniques are used to show that the Coxeter
matrix $\Phi$ is periodic if $\Bform$ is positive and weakly periodic
if $\Bform$ is non-negative. It is much harder to establish
implications in the other direction. As already noted in~\cite{Sato05},
even if $\Phi$ is periodic, $\Bform$ may be indefinite, so additional
constraints are needed.

An alternative definition of the Coxeter matrix is as a certain product
of reflections defined by a generalized Cartan
matrix~\cite{ACampo76,Ringel94}, whereas the definition given above is
$-C^{-1} C^t$ where $C$ is the matrix of the bilinear form.

We claim similarly to~\cite{Coleman89}, and give a simple, direct
proof, that for any square matrix $A$ with $2$ on its main diagonal,
the product of the $n$ reflections it defines can be expressed as
$-A_{+}^{-1}A_{-}^t$ where $A_{+}$, $A_{-}$ are closely associated with
the upper and lower triangular parts of $A$, see
Section~\ref{sec:matrices}. This claim can be generalized to products
in arbitrary order, and no other conditions on $A$, such as being
generalized Cartan, bipartite~\cite{ACampo76} or
symmetric~\cite{Howlett82}, are needed. In particular, when $\Bform$ is
triangular, $\Phi$ can be written as a product of the reflections
defined by the symmetrization of $\Bform$.

Further connections between periodicity and non-negativity are achieved
when we restrict ourselves to pairs $(V, \Bform)$ for which there
exists a finite dimensional $k$-algebra $\gL$ over an algebraically
closed field $k$, having finite global dimension, such that $V \cong
K_0(\Mod \gL)$ and $\Bform$ coincides, under that isomorphism, with the
Euler form $\Bform_{\gL}$ of $\gL$. Here $\Mod \gL$ denotes the
category of finite dimensional right $\gL$-modules. Since $\gldim \gL <
\infty$, the form $\Bform_{\gL}$ is non-degenerate, hence its Coxeter
transformation $\Phi_{\gL}$ is well-defined and coincides with the
image in $K_0(\Mod \gL)$ of the Auslander-Reiten translation on the
bounded derived category $\cD^b(\Mod \gL)$.

In Section~\ref{sec:algebras} of the paper, we show that if $\gL$ is
\emph{piecewise hereditary}, i.e. its bounded derived category
$\cD^b(\Mod \gL)$ is equivalent as a triangulated category to
$\cD^b(\cH)$ for a hereditary abelian category $\cH$, then the
periodicity of $\Phi_\gL$ implies the non-negativity of $\Bform_\gL$.

In that Section, we also show that when $\gL$ is an incidence algebra
of a poset $X$, the Euler form $\Bform_{\gL}$ and its Coxeter
transformation $\Phi_{\gL}$ can be explicitly described in terms of the
combinatorics of $X$.

Previously, \cite{dlPena03} claimed that the condition of $\Bform$
being \emph{triangular}, that is, its matrix with respect to some basis
of $V$ is upper triangular with ones on the main diagonal, is enough
for the periodicity of $\Phi$ to imply the non-negativity of $\Bform$.
We find however examples of incidence algebras of posets negating this
claim, see Section~\ref{sec:examples}.

\section{Coxeter transformations of bilinear forms}
\label{sec:matrices}

\subsection{The definition of the Coxeter matrix}
Let $V$ be a free abelian group of finite rank and let $\Bform : V
\times V \to \bZ$ be a non-degenerate $\bZ$-bilinear form on $V$.
Recall that $\Bform$ is \emph{positive} if $\bform{v}{v} > 0$ for all
$0 \neq v \in V$, \emph{non-negative} if $\bform{v}{v} \geq 0$ for all
$v \in V$ and \emph{indefinite} otherwise. The \emph{Coxeter
transformation} $\Phi : V \to V$ corresponding to $\Bform$ is defined
via the equation $\bform{v}{w} = - \bform{w}{\Phi v}$ for all $v, w \in
V$~\cite{Lenzing96}.

We consider the elements of $\bZ^n$ as column vectors, and denote by
$M^{t}$ the transpose of a matrix $M$. Let $\{e_i\}_{i=1}^{n}$ be the
standard basis of $\bZ^n$. By choosing a $\bZ$-basis $v_1,\dots,v_n$ of
$V$, we may identify $V$ with $\bZ^n$ and $\Bform$ with the form
$\Bform_C$ defined by $\bform{x}{y}_C = x^t C y$ where $C \in
\GL_n(\bZ)$ is the matrix whose entries are $C_{ij} = \bform{v_i}{v_j}$
for $1 \leq i,j \leq n$. In other words, $\bform{v_i}{v_j} =
\bform{e_i}{e_j}_C$. Under this identification, the matrix of $\Phi$ is
$-C^{-1}C^t$, hence we define the \emph{Coxeter matrix} $\Phi_C$ of a
matrix $C \in \GL_n(\bZ)$ to be $\Phi_C = -C^{-1}C^t$.

Note that $v = -C^{-1}C^t v$ if and only if $(C + C^t)v = 0$, hence the
geometric multiplicity of the eigenvalue $1$ in $\Phi_C$ equals the
dimension of the radical of the symmetrized bilinear form $C+C^t$.

\begin{defn}
A matrix $\Phi \in \GL_n(\bZ)$ is \emph{periodic} if $\Phi^m=I$ for
some $m \geq 1$. $\Phi$ is \emph{weakly periodic} if for some $m
\geq 1$, $\Phi^m-I$ is nilpotent.
\end{defn}

\begin{defn}
A matrix $C \in \GL_n(\bZ)$ is \emph{unitriangular} if $C$ is upper
triangular and $C_{ii}=1$ for $1 \leq i \leq n$.
\end{defn}

Relations between the positivity of the bilinear form $\Bform_C$ and
the periodicity of $\Phi_C$ have been studied in~\cite{dlPena03,Sato05}
and are summarized as follows:

\begin{thm}
Let $C \in \GL_n(\bZ)$. Then:
\begin{enumerate}
\item \cite[(2.8)]{Sato05} $\Phi_C$ is periodic if $\Bform_C$ is
positive.
\item \cite[(3.4)]{Sato05} $\Phi_C$ is weakly periodic if
$\Bform_C$ is non-negative.
\end{enumerate}
\end{thm}

However, \cite[(3.8)]{Sato05} is an example of a matrix whose Coxeter
matrix is periodic but the corresponding bilinear form is indefinite.

\subsection{Alternative definition as a product of reflections}
Following \cite{ACampo76,BLM89,Ringel94}, we review an alternative
definition of the Coxeter matrix as a product of reflections.

Let $A$ be an $n \times n$ matrix with integer entries satisfying
\begin{align}
\tag{A1} \label{e:A1} A_{ii} &= 2 \qquad 1\leq i \leq n \\
\tag{A2} \label{e:A2} A_{ij} &= 0 \qquad \text{if and only if }
A_{ij} = 0, 1 \leq i, j \leq n
\end{align}
The \emph{primitive graph} of $A$ (cf.~\cite{BLM89}) is an
undirected graph with $n$ vertices, where two vertices $i \neq j$
are connected by an edge if $A_{ij} \neq 0$. The matrix $A$ is
\emph{indecomposable} if its primitive graph is connected.

Define \emph{reflections} $r_1, \dots, r_n$ by
\begin{equation}
\label{e:reflect} r_i(e_j) = e_j - A_{ij}e_i \qquad 1 \leq j \leq n
\end{equation}
In other words, $r_i$ is the matrix obtained from the identity matrix
by subtracting the $i$-th row of $A$. Denote by $I$ the $n \times n$
identity matrix.

\begin{lem}
\label{l:ri} Let $A$ be a matrix satisfying \eqref{e:A1}.
\begin{enumerate}
\renewcommand{\theenumi}{\alph{enumi}}
\item
$r_i^2 = I$ for $1 \leq i \leq n$.
\item
If $A$ satisfies also \eqref{e:A2}, then $r_i r_j = r_j r_i$ for any
two non-adjacent vertices $i, j$ on the primitive graph of $A$.
\end{enumerate}
\end{lem}
\begin{proof}
Since $A_{ii} = 2$, we have $r_i(e_i) = -e_i$, thus
\[
r_i^2(e_t) = r_i(e_t - A_{it}e_i) = e_t - A_{it}e_i - A_{it}r_i(e_i)
= e_t
\]
for all $1 \leq t \leq n$, and the first assertion is proved.

If $A_{ij}=0$ then $r_i(e_j) = e_j$. The assumptions on $A$ imply
that if $i, j$ are not adjacent, then $r_i(e_j) = e_j$ and $r_j(e_i)
= e_i$. Therefore, if $1 \leq t \leq n$,
\[
r_i r_j(e_t) = r_i(e_t - A_{jt}e_j) = e_t - A_{it}e_i - A_{jt}e_j
\]
is symmetric in $i$ and $j$, hence $r_i r_j = r_j r_i$.
\end{proof}

Consider the following two additional properties:
\begin{align}
\tag{A3} \label{e:A3} &A_{ij} \leq 0 \qquad \text{for all $i \neq j$} \\
\tag{A4} \label{e:A4} &\text{The primitive graph of $A$ is
bipartite}
\end{align}
\begin{defn}
A matrix $A$ is a \emph{generalized Cartan matrix} if it satisfies
\eqref{e:A1}, \eqref{e:A2} and \eqref{e:A3}. A matrix $A$ is
\emph{bipartite} if it satisfies \eqref{e:A1}, \eqref{e:A2} and
\eqref{e:A4}.
\end{defn}

For a generalized Cartan matrix $A$ and a permutation $\pi$ of
$\{1,2,\dots,n\}$, a \emph{Coxeter transformation} is defined
in~\cite{Ringel94} by $\Phi(A,\pi) = r_{\pi(1)}r_{\pi(2)} \cdots
r_{\pi(n)}$. For a bipartite matrix $A$, let $\Sigma_1 \amalg \Sigma_2$
be a corresponding partition of $\{1,2,\dots,n\}$ and consider $R_A =
R_1 R_2$ where $R_k = \prod_{i \in \Sigma_k} r_i$, $k=1,2$,
see~\cite{BLM89}. Note that by Lemma~\ref{l:ri}, the matrices $R_k$ do
not depend on the order of reflections within each product. Note also
that $R_A$ equals $\Phi(A,\pi)$ for a suitable $\pi$.


Recall that the \emph{spectrum} of a square matrix $\Phi$ with complex
entries, denoted $\spec(\Phi)$, is the set of (complex) roots of the
characteristic polynomial of $\Phi$. Let $\rho(\Phi) = \max \left\{
|\lambda| \,:\, \lambda \in \spec(\Phi) \right\}$ be the \emph{spectral
radius} of $\Phi$.

We recall two results on the spectrum of Coxeter transformations
corresponding to generalized Cartan and bipartite matrices.


\begin{thm}[\cite{Ringel94}] \label{t:gcm}
Let $A$ be an indecomposable generalized Cartan matrix, $\pi \in S_n$.
If $A$ is not of finite or affine type, then $\rho(\Phi(A,\pi)) > 1$.
\end{thm}

\begin{thm}[\protect{\cite[p.~63]{ACampo76},\cite[p.~344]{BLM89}}]
\label{t:bipartite} Let $A$ be a bipartite matrix.
\begin{enumerate}
\renewcommand{\theenumi}{\alph{enumi}}
\item
$ \lambda^2 \in \spec(R_A) \text{ if and only if } \lambda + 2 +
\lambda^{-1} \in \spec(A)$.
\item
If $A$ is also symmetric, then $\spec(R_A) \subset S^1 \cup \bR$.
\end{enumerate}
\end{thm}

\subsection{Linking the two definitions}

Let $R$ be any commutative ring with $1$ and let $e_1,\dots,e_n$ be a
basis of a free $R$-module of rank $n$. Let $A$ be an $n \times n$
matrix with entries in $R$ satisfying~\eqref{e:A1} (where $2$ means
$1+1$), and define the reflections $r_1, \dots, r_n$ as
in~\eqref{e:reflect}. When we want to stress the dependence of the
reflections on $A$, we shall use the notation $r^A_1,\dots,r^A_n$.

\begin{lem}
\label{l:prod} Let $1 \leq s \leq n$. Then for every $1 \leq t \leq n$,
\[
(r_1 \cdots r_s)(e_t) = e_t + \sum_{k=1}^{s} (-1)^k \sum_{1 \leq i_1 <
i_2 < \dots < i_k \leq s} A_{i_1 i_2} \cdots A_{i_{k-1} i_k} A_{i_k t}
e_{i_1}
\]
\end{lem}
\begin{proof}
By induction on $s$, the case $s=1$ being just the definition of $r_1$,
and for the induction step, expand $r_{s+1}(e_t)$ as $e_t - A_{s+1,t}
e_{s+1}$ and use the hypothesis for $s$.
\[
\begin{split}
&(r_1 \cdots r_s r_{s+1})(e_t) = (r_1 \cdots r_s)(e_t) - A_{s+1,t}
(r_1 \cdots r_s)(e_{s+1}) \\
&= e_t - A_{s+1,t} e_{s+1} + \sum_{k=1}^{s} (-1)^k \sum_{1 \leq i_1 <
i_2 < \dots < i_k \leq s} A_{i_1 i_2} \cdots A_{i_{k-1} i_k} A_{i_k t}
e_{i_1}
\\ &\quad
+ \sum_{k=1}^{s} (-1)^{k+1} \sum_{1 \leq i_1 < i_2 < \dots < i_k
\leq s} A_{i_1 i_2} \cdots A_{i_{k-1} i_k} A_{i_k,s+1} A_{s+1,t} e_{i_1}
\\
&= e_t - \sum_{k=1}^{s+1} (-1)^k \sum_{1 \leq i_1 < i_2 < \dots < i_k
\leq {s+1}} A_{i_1 i_2} \cdots A_{i_{k-1} i_k} A_{i_k t} e_{i_1}
\end{split}
\]
\end{proof}

Define two $n \times n$ matrices $A_+$ and $A_{-}$ by
\begin{align*}
(A_{+})_{ij} = \begin{cases} A_{ij} & i < j \\ 1 & i = j \\ 0 & i > j
\end{cases} & & (A_{-})_{ij} = \begin{cases} A_{ji} & i < j \\ 1 & i =
j \\ 0 & i > j \end{cases}
\end{align*}

Then $A = A_+ + A_{-}^t$, and one can think of $A_+$, $A_{-}$ as the
upper and lower triangular parts of $A$. The matrices $A_{+}$ and
$A_{-}$ are invertible since $A_{+}-I$ and $A_{-}-I$ are nilpotent.
Note that $A$ is symmetric if and only if $A_+ = A_{-}$

\begin{thm}
\label{t:prodref} If $A$ satisfies~\eqref{e:A1}, then $r^A_1 r^A_2
\cdots r^A_n = -A_+^{-1} A_{-}^t$.
\end{thm}
\begin{proof}
By Lemma~\ref{l:prod} with $s=n$,
\[
(r_1 \cdots r_n)(e_t) = e_t + \sum_{k=1}^{n} (-1)^k \sum_{1 \leq i_1 <
i_2 < \dots < i_k \leq n} A_{i_1 i_2} \cdots A_{i_{k-1} i_k} A_{i_k t}
e_{i_1}
\]

This can be written in matrix form, using the definition of $A_{+}$, as
follows:
\begin{align*}
r_1 \cdots r_n &= I + \sum_{k=1}^{n} (-1)^k
(A_{+}-I)^{k-1} A \\
&= I - (I - (A_{+}-I) + (A_{+}-I)^2 - \dots) A \\ &= I - A_{+}^{-1}
(A_{+} + A_{-}^t) = -A_{+}^{-1} A_{-}^t
\end{align*}
\end{proof}

\begin{rem}
Theorem~\ref{t:prodref} is still true when we drop the
condition~\eqref{e:A1} and slightly change the definition of $A_{-}$,
by $(A_{-})_{ii} = A_{ii}-1$ for $1 \leq i \leq n$. However, in that
case the matrices $r_i$ are no longer reflections.
\end{rem}

Theorem~\ref{t:prodref} provides a link between the definition of the
Coxeter matrix as a specific automorphism of the bilinear form and its
definition as a product of $n$ reflections, as shown by the following
corollary.

\begin{cor} 
\label{c:prodref} Let $C \in \GL_n(\bZ)$ be a unitriangular matrix.
Then $\Phi_C = \Phi(C+C^t, id)$, that is, $\Phi_C = r^A_1 r^A_2 \cdots
r^A_n$ for $A = C + C^t$.
\end{cor}

In fact, this corollary is proved in~\cite{Howlett82} for the case
where $\Phi_C$ is a Coxeter element in an arbitrary Coxeter group of
finite rank represented as a group of linear transformations on a real
inner product space, so that the Cartan matrix $A$ is symmetric.

\begin{proof}
Apply Theorem~\ref{t:prodref} for the matrix $A = C + C^t$, which
satisfies \eqref{e:A1}, \eqref{e:A2} and $A_{+} = A_{-} = C$.
\end{proof}

Denote by $S_n$ the group of permutations on $\{1,2,\dots,n\}$ and let
$\pi \in S_n$. One could deduce a generalized version of
Theorem~\ref{t:prodref} for the product of the $n$ reflections in an
arbitrary order by proving an analogue of Lemma~\ref{l:prod} for
arbitrary $\pi$. Instead, we will derive the generalized version from
the original one using permutation matrices.

Define the permutation matrix $P_\pi$ by $P_\pi(e_i) = e_{\pi(i)}$ for
all $1 \leq i \leq n$. Note that $P_\pi^{-1} = P_\pi^t$. Given a matrix
$A$, let $A_\pi$ denote the matrix $P_\pi^{-1} A P_\pi$, so that
$(A_\pi)_{ij} = A_{\pi(i) \pi(j)}$. Obviously, if $A$
satisfies~\eqref{e:A1}, so does $A_\pi$.

\begin{lem}
\label{l:rApi} Let $1 \leq i \leq n$. Then $r^{A_\pi}_i = P_\pi^{-1}
r^A_{\pi(i)} P_\pi$.
\end{lem}
\begin{proof}
For all $1 \leq t \leq n$,
\[
\left(P_\pi^{-1} r^A_{\pi(i)} P_\pi \right)(e_t) = P_\pi^{-1}
\left(e_{\pi(t)} - A_{\pi(i) \pi(t)} e_{\pi(i)} \right) = e_t -
A_{\pi(i)\pi(t)} e_i = r^{A_\pi}_i(e_t)
\]
\end{proof}

Define two $n \times n$ matrices $A_{\pi,+}$ and $A_{\pi,-}$ by
\begin{align*}
(A_{\pi,+})_{ij} = \begin{cases} A_{ij} & \pi^{-1}(i) < \pi^{-1}(j) \\
1 & i=j \\
0 & \text{otherwise}
\end{cases}
& &
(A_{\pi,-})_{ij} = \begin{cases} A_{ji} & \pi^{-1}(i) < \pi^{-1}(j) \\
1 & i=j \\
0 & \text{otherwise}
\end{cases}
\end{align*}
Direct calculation shows that $A_{\pi,+} = P_\pi (A_\pi)_{+}
P_\pi^{-1}$, $A_{\pi,-} = P_\pi (A_\pi)_{-} P_\pi^t$ and $A = A_{\pi,+}
+ A_{\pi,-}^t$.

\begin{cor}
Let $A$ satisfy~\eqref{e:A1} and let $\pi \in S_n$. Then
\[
r^A_{\pi(1)} r^A_{\pi(2)} \dots r^A_{\pi(n)} = - A_{\pi,+}^{-1}
A_{\pi,-}^t
\]
\end{cor}
\begin{proof}
By Lemma~\ref{l:rApi} and Theorem~\ref{t:prodref} applied for $A_\pi$,
\begin{align*}
r^A_{\pi(1)} r^A_{\pi(2)} \dots r^A_{\pi(n)} &= P_\pi \left(
r^{A_\pi}_1 r^{A_\pi}_2 \dots r^{A_\pi}_n \right) P_\pi^{-1} = - P_\pi
(A_\pi)_{+}^{-1} (A_\pi)_{-}^t P_\pi^t \\
&= - \left( P_\pi (A_\pi)_{+}^{-1} P_\pi^{-1} \right) \left( P_\pi
(A_\pi)_{-}^t P_\pi^t \right) = -A_{\pi,+}^{-1} A_{\pi,-}^t
\end{align*}
\end{proof}

\section{Periodicity and non-negativity for piecewise hereditary algebras and posets}
\label{sec:algebras}

Let $k$ be a field, and let $\cA$ be an abelian $k$-category of finite
global dimension with finite dimensional $\Ext$-spaces. Denote by
$\cD^b(\cA)$ its bounded derived category and by $K_0(\cA)$ its
Grothendieck group. The expression
\[
\bform{X}{Y}_{\cA} = \sum_{i \in \bZ} (-1)^i \dim_k
\Hom_{\cD^b(\cA)}(X, Y[i])
\]
is well-defined for $X, Y \in \cD^b(\cA)$ and induces a $\bZ$-bilinear
form on $K_0(\cA)$, known as the \emph{Euler form}. When $\Bform_{\cA}$
is non-degenerate, the unique transformation $\Phi_\cA : K_0(\cA) \to
K_0(\cA)$ satisfying $\bform{x}{y}_{\cA} = -\bform{y}{\Phi_\cA
x}_{\cA}$ for all $x, y \in K_0(\cA)$ is called the \emph{Coxeter
transformation} of $\cA$. For more details we refer the reader
to~\cite{Lenzing99}.

Two such abelian $k$-categories $\cA$ and $\cB$ are said to be
\emph{derived equivalent} if there exists a triangulated equivalence $F
: \cD^b(\cA) \simeq \cD^b(\cB)$. In this case, the forms $\Bform_{\cA}$
and $\Bform_{\cB}$ are equivalent over $\bZ$, hence the positivity
properties of the Euler form and the periodicity properties of the
Coxeter transformation are invariants of derived equivalence.

Let $\gL$ be a finite dimensional algebra of finite global dimension
over an algebraically closed field $k$, and consider the $k$-category
$\Mod \gL$ of finitely generated right modules over $\gL$. Denote by
$\dL$ its bounded derived category, by $K_0(\gL)$ its Grothendieck
group and by $\Bform_{\gL}$ the Euler form. Then $K_0(\gL)$ is free of
finite rank, with a $\bZ$-basis consisting of the representatives of
the isomorphism classes of simple modules in $\Mod \gL$. The form
$\Bform_\gL$ is non-degenerate, and its Coxeter transformation
$\Phi_\gL$ coincides with the linear map on $K_0(\gL)$ induced by the
Auslander-Reiten translation on $\dL$. For more details
see~\cite[(III.1)]{Happel88}, \cite[(2.4)]{Ringel84}
or~\cite{Lenzing99}.

\subsection{Path algebras of quivers without oriented cycles}

The first example of algebras $\gL$ for which the connection between
the positivity of $\Bform_\gL$ and the periodicity of $\Phi_\gL$ is
completely understood is the class of path algebras of quivers without
oriented cycles, or more generally hereditary algebras,
see~\cite[Theorem~18.5]{Lenzing99}. We briefly review the main results.

A (finite) \emph{quiver} $Q$ is a directed graph with a finite number
of vertices and edges. The \emph{underlying graph} of $Q$ is the
undirected graph obtained from $Q$ by forgetting the orientations of
the edges.  An \emph{oriented cycle} is a nontrivial path in $Q$
starting and ending at the same vertex. The \emph{path algebra} $kQ$ is
the algebra over $k$ having as a $k$-basis the set of all (oriented)
paths in $Q$; the product of two paths is their composition, if
defined, and zero otherwise.

When $Q$ has no oriented cycles, the path algebra $kQ$ is hereditary
and finite-dimensional. Denote by $\Bform_Q$ its Euler form and by
$\Phi_Q$ its Coxeter transformation. The matrix of $\Bform_Q$ with
respect to the basis of simple modules is unitriangular, and its
symmetrization is generalized Cartan. The relations between the
periodicity of $\Phi_Q$ and the positivity of $\Bform_Q$ are summarized
in the following well-known proposition,
see~\cite{ACampo76,BGP73,Bourbaki02,Ringel94}
and~\cite[(1.2)]{Ringel84}.

\begin{prop} \label{p:quivercase}
Let $Q$ be a connected quiver without oriented cycles. Then:
\begin{enumerate}
\renewcommand{\theenumi}{\alph{enumi}}
\item
$\Phi_Q$ is periodic if and only if $\Bform_Q$ is positive,
equivalently the underlying graph of $Q$ is a Dynkin diagram of type
$A$, $D$ or $E$.

\item
$\Phi_Q$ is weakly periodic if and only if $\Bform_Q$ is non-negative,
equivalently the underlying graph of $Q$ is a Dynkin diagram or an
extended Dynkin diagram of type $\tilde{A}$, $\tilde{D}$ or
$\tilde{E}$.
\end{enumerate}
\end{prop}

\subsection{Canonical algebras}

Another interesting class of algebras for which the connection between
non-negativity and periodicity is established are the canonical
algebras, introduced in~\cite{Ringel84}.

The Grothendieck group and the Euler form of canonical algebras were
thoroughly studied in~\cite{Lenzing96}. If $\gL$ is canonical of type
$(\bp,\bl)$ where $\bp=(p_1,\dots,p_t)$ and $\bl=(\lambda_3, \dots,
\lambda_t)$ is a sequence of pairwise distinct elements of $k \setminus
\{0\}$, then the rank of $K_0(\gL)$ is $\sum_{i=1}^t p_i - (t-2)$ and
the characteristic polynomial of the Coxeter transformation
$\Phi_{\gL}$ equals $(T-1)^2 \prod_{i=1}^t
\frac{T^{p_i}-1}{T-1}$~\cite[Prop.~7.8)]{Lenzing96}. In particular,
$\rho(\Phi_{\gL})=1$ and the eigenvalues of $\Phi_{\gL}$ are roots of
unity, hence $\Phi_{\gL}$ is weakly periodic.

The following proposition follows from~\cite[Prop.~10.3]{Lenzing96},
see also~\cite{LenzingdlPena97}.

\begin{prop}
\label{p:canonicalcase} Let $\gL$ be a canonical algebra of type
$(\bp,\bl)$. If $\Phi_{\gL}$ is periodic then $\bp$ is one of
$(2,3,6)$, $(2,4,4)$, $(3,3,3)$ or $(2,2,2,2)$. In any of these cases,
$\Bform_{\gL}$ is non-negative.
\end{prop}

\subsection{Extending to piecewise hereditary algebras}

We extend the results of the previous sections to the class of all
piecewise hereditary algebras.

\begin{defn}
An algebra $\gL$ over $k$ is \emph{piecewise hereditary} if there exist
a hereditary abelian category $\cH$ and a triangulated equivalence $\dL
\simeq \dH$.
\end{defn}

\begin{thm}
Let $k$ be an algebraically closed field and let $\gL$ be a finite
dimensional piecewise hereditary $k$-algebra. If $\Phi_{\gL}$ is
periodic, then $\Bform_{\gL}$ is non-negative.
\end{thm}

\begin{proof}
By definition, there exists a hereditary category $\cH$ and an
equivalence of triangulated categories $F: \dL \simeq \dH$. By the
invariance under derived equivalence, it is enough to prove the theorem
for $\Phi_\cH$ and $\Bform_{\cH}$. Moreover, we can assume that $\cH$
is connected.

Now $\cH$ is an $\Ext$-finite $k$-category and $F(\gL_{\gL})$ is a
tilting complex in $\dH$, so by~\cite[Theorem~1.7]{HappelReiten98},
$\cH$ admits a \emph{tilting object}, that is, an object $T$ with
$\Ext^1_{\cH}(T, T) = 0$ such that for any object $X$ of $\cH$, the
condition $\Hom_{\cH}(T,X) = 0 = \Ext^1_{\cH}(T,X)$ implies that $X =
0$.

By the classification of hereditary connected $\Ext$-finite
$k$-categories with tilting object up to derived equivalence over an
algebraically closed field~\cite{Happel01}, $\cH$ is derived equivalent
to $\Mod H$ for a finite dimensional hereditary algebra $H$ or to $\Mod
\gL$ for a canonical algebra $\gL$. Again by invariance under derived
equivalence we may assume that $\cH = \Mod H$ or $\cH = \Mod \gL$.

For $\cH = \Mod H$, we can replace $H$ by a path algebra of a finite
connected quiver without oriented cycles, and then use
Proposition~\ref{p:quivercase}. For $\cH = \Mod \gL$, the result
follows from Proposition~\ref{p:canonicalcase}.
\end{proof}

\subsection{Incidence algebras of posets}

Let $X$ be a finite partially ordered set (\emph{poset}) and let $k$ be
a field. The \emph{incidence algebra} $kX$ is the $k$-algebra spanned
by elements $e_{xy}$ for the pairs $x \leq y$ in $X$, with
multiplication defined by $e_{xy} e_{zw} = \delta_{yz} e_{xw}$. Finite
dimensional right modules over $kX$ can be identified with commutative
diagrams of finite dimensional $k$-vector spaces over the \emph{Hasse
diagram} of $X$ which is the directed graph whose vertices are the
points of $X$, with an arrow from $x$ to $y$ if $x < y$ and there is no
$z \in X$ with $x < z < y$.

We recollect the basic facts on the Euler form of posets and refer the
reader to~\cite{Ladkani07} for details. The algebra $kX$ is of finite
global dimension, hence its Euler form, denoted $\Bform_X$, is
well-defined and non-degenerate. Denote by $C_X$, $\Phi_X$ the matrices
of $\Bform_X$ and its Coxeter transformation with respect to the basis
of simple $kX$-modules.

The \emph{incidence matrix} of $X$, denoted $\one_X$, is the $X \times
X$ matrix defined by
\[
(\one_X)_{xy} = \begin{cases} 1 & x \leq y \\ 0 & \text{otherwise}
\end{cases}
\]
By extending the partial order on $X$ to a linear order, we can always
arrange the elements of $X$ such that the incidence matrix is
unitriangular. In particular, $\one_X$ is invertible over $\bZ$. Recall
that the \emph{M\"obius function} $\mu_X : X \times X \to \bZ$ is
defined by $\mu_X(x,y) = (\one_X^{-1})_{xy}$.

\begin{lem}[\protect{\cite[Prop.~3.11]{Ladkani07}}]
$C_X = \one_X^{-1}$.
\end{lem}

\begin{lem}
\label{l:coxposet} Let $x, y \in X$. Then $(\Phi_X)_{xy} = -\sum_{z
\,:\, z \geq x} \mu_X(y,z)$.
\end{lem}
\begin{proof}
Since $\Phi_X = -C_X^{-1} C_X^t = -\one_X \one_X^{-t}$,
\[
(\Phi_X)_{xy} = - \sum_{z \in X} (\one_X)_{xz} (\one_X^{-1})_{yz} =
- \sum_{z \,:\, z \geq x} \mu_X(y, z)
\]
\end{proof}

When the Hasse diagram of $X$ has the property that any two vertices
$x, y$ are connected by at most one directed path, the M\"obius
function takes a very simple form, namely
\[
\mu_X(x,y) = \begin{cases} 1 & y=x \\ -1 & \text{$x \to y$ is an
edge in the Hasse diagram} \\ 0 & \text{otherwise} \end{cases}
\]
In this case, Lemma~\ref{l:coxposet} coincides with Proposition~3.1
of~\cite{Boldt95}, taking the Hasse diagram as the quiver.

\begin{lem}
If $X$ and $Y$ are posets, then
\[
C_{X \times Y} = C_X \otimes C_Y \qquad \Phi_{X \times Y} = -\Phi_X
\otimes \Phi_Y
\]
\end{lem}
\begin{proof}
Observe that $\one_{X \times Y} = \one_X \otimes \one_Y$.
\end{proof}

\begin{cor}
\label{c:coxprod} Let $X$, $Y$ be posets with periodic Coxeter
matrices. Then $X \times Y$ has also periodic Coxeter matrix.
\end{cor}

Since non-negativity of forms is not preserved under tensor products,
Corollary~\ref{c:coxprod} can be used to construct posets with periodic
Coxeter matrix but with indefinite Euler form, see
Example~\ref{ex:prod}.

\section{Examples}
\label{sec:examples}
For a poset $X$, let $C_X$, $\Phi_X$ be as in the previous section. In
particular we may assume that $C_X$ is unitriangular. The
symmetrization $A_X = C_X + C_X^t$ satisfies~\eqref{e:A1}
and~\eqref{e:A2}, but in general it is not bipartite nor generalized
Cartan.

\subsection{Spectral properties of $\Phi_X$}

\begin{exmp}
\label{ex:specM} \emph{The spectrum of $\Phi_X$ does not determine
that of $A_X$ (Compare with Theorem~\ref{t:bipartite}a).}

The four posets whose Hasse diagrams are depicted in
Figure~\ref{fig:specM} are derived equivalent (as they are all
piecewise hereditary of type $D_5$), hence their Coxeter matrices are
similar and have the same spectrum, namely the roots of the
characteristic polynomial $x^5+x^4+x+1$. However, the spectra of the
corresponding symmetrized forms are different. Figure~\ref{fig:specM}
also shows for each poset $X$ the characteristic polynomial of the
matrix of its symmetrized form.

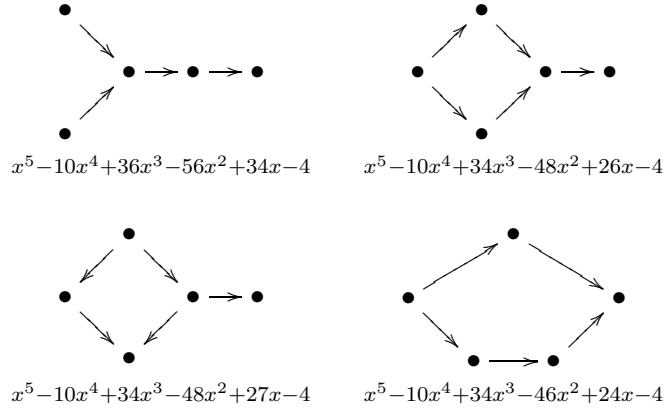
\begin{figure}
\[
\begin{array}{cccc}
\xymatrix@=1pc{
{\bullet} \ar[dr] \\
& {\bullet} \ar[r] & {\bullet} \ar[r] & {\bullet} \\
{\bullet} \ar[ur]
} & & \xymatrix@=1pc{
& {\bullet} \ar[dr] \\
{\bullet} \ar[dr] \ar[ur] & & {\bullet} \ar[r] & {\bullet} \\
& {\bullet} \ar[ur]
}
\\
\scriptstyle{x^5-10x^4+36x^3-56x^2+34x-4} & &
\scriptstyle{x^5-10x^4+34x^3-48x^2+26x-4} \\ \\
\xymatrix@=1pc{
& {\bullet} \ar[dl] \ar[dr] \\
{\bullet} \ar[dr] & & {\bullet} \ar[dl] \ar[r] & {\bullet} \\
& {\bullet}
} & & \xymatrix@=0.25pc{
&&& {\bullet} \ar[ddrrr] \\ \\
{\bullet} \ar[uurrr] \ar[ddrr] &&&&&& {\bullet} \\ \\
&& {\bullet} \ar[rr] && {\bullet} \ar[uurr]
}
\\
\scriptstyle{x^5-10x^4+34x^3-48x^2+27x-4} & &
\scriptstyle{x^5-10x^4+34x^3-46x^2+24x-4}
\end{array}
\]
\caption{Derived equivalent posets with different spectra of the
corresponding symmetrized bilinear forms.} \label{fig:specM}
\end{figure}
\end{exmp}

\begin{exmp}
\label{ex:spec} \emph{A poset $X$ with $\spec{\Phi_X} \not \subseteq
S^1 \cup \bR$ (Compare with Theorem~\ref{t:bipartite}b).}

Let $X$ be the following poset.
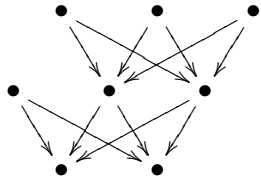
\[
\xymatrix@=0.5pc{ & {\bullet} \ar[ddr] \ar[ddrrr] &&
{\bullet} \ar[ddl] \ar[ddr] && {\bullet} \ar[ddl] \ar[ddlll] \\
\\
{\bullet} \ar[ddr] \ar[ddrrr] &&
{\bullet} \ar[ddl] \ar[ddr] && {\bullet} \ar[ddl] \ar[ddlll] \\
\\
& {\bullet} && {\bullet} }
\]
The characteristic polynomial of $\Phi_X$ is
$(x+1)^4(x^4-2x^3+6x^2-2x+1)$, whose roots, besides $-1$, are $z,
\bar{z}, z^{-1}, \bar{z}^{-1}$ with $\Re{z} =
\frac{1+\sqrt{2\sqrt{3}-3}}{2}$ and $|z|^2 =
\frac{1+\sqrt{2\sqrt{3}-3}}{2-\sqrt{3}} - 1$. These four roots are
neither real nor on the unit circle.

An example of similar spectral behavior for path algebra of a quiver is
given in~\cite[Example~18.1]{Lenzing99}.

Note that for all posets $X$ with 7 elements or less, $\spec(\Phi_X)
\subseteq S^1 \cup \bR $. This was verified using the
database~\cite{PosetAtlas} and the \textsc{Magma} software package.
\end{exmp}

\subsection{Counterexamples to~\cite[Prop.~1.2]{dlPena03}}

We give two examples of posets showing that in general, for triangular
algebras, the periodicity of the Coxeter transformation (and even of
the Auslander-Reiten translation up to a shift) does not imply the
non-negativity of the Euler form.

\begin{exmp}
\label{ex:periodic}

Consider the poset $X$ with the following Hasse diagram.
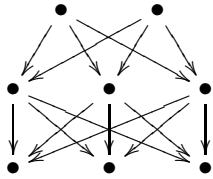
\[
\xymatrix@=0.5pc{ & {\bullet} \ar[ddl] \ar[ddr] \ar[ddrrr] &&
{\bullet} \ar[ddlll] \ar[ddl] \ar[ddr] \\ \\
{\bullet} \ar[dd] \ar[ddrr] \ar[ddrrrr] && {\bullet} \ar[ddll]
\ar[dd] \ar[ddrr] &&
{\bullet} \ar[ddllll] \ar[ddll] \ar[dd] \\ \\
{\bullet} && {\bullet} && {\bullet} }
\]
Then $\Phi_X^6=I$ but $v^t C_X v = -1$ for $v= \begin{pmatrix}
1&1&1&1&1&0&0&0 \end{pmatrix}^t$ (the vertices are ordered in layers
from top to bottom).
\end{exmp}

\begin{exmp}
\label{ex:prod}

Let $X = A_3 \times D_4$ with the following orientations:
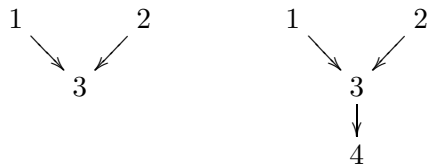
\[
\xymatrix@=1pc{
{1} \ar[dr] & & {2} \ar[dl] \\
& {3} } \qquad \qquad \xymatrix@=1pc{
{1} \ar[dr] & & {2} \ar[dl] \\
& {3} \ar[d] \\
& {4} }
\]
The Hasse diagram of $X$ is given by
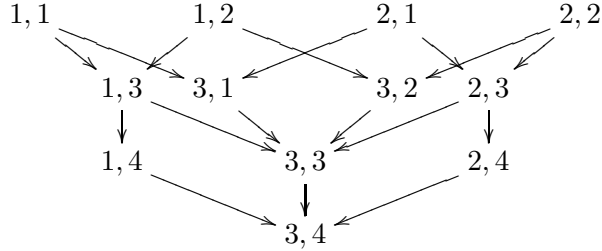
\[
\xymatrix@=1pc{ {1,1} \ar[dr] \ar[drr] & & {1,2} \ar[dl] \ar[drr] &
&
{2,1} \ar[dr] \ar[dll] & & {2,2} \ar[dl] \ar[dll] \\
& {1,3} \ar[d] \ar[drr] & {3,1} \ar[dr] & &
{3,2} \ar[dl] & {2,3} \ar[d] \ar[dll]\\
& {1,4} \ar[drr] & & {3,3} \ar[d] & & {2,4} \ar[dll] \\
& & & {3,4} }
\]
so that $X$ contains the following wild quiver as a subposet.
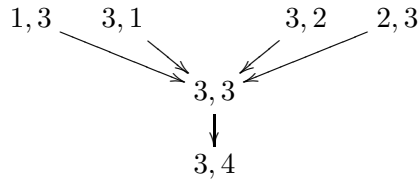
\[
\xymatrix@=1pc{ {1,3} \ar[drr] & {3,1} \ar[dr] & &
{3,2} \ar[dl] & {2,3} \ar[dll]\\
& & {3,3} \ar[d] \\
& & {3,4} }
\]
It follows~\cite{Loupias75} that $kX$ is not of finite representation
type, hence by~\cite[Theorem~6]{Yuzvinsky81} the form $\Bform_X$ is not
weakly positive, that is, there exists a vector $v \neq 0$ with
non-negative coordinates such that $\bform{v}{v}_X \leq 0$.

Moreover, we can exhibit a non-negative vector $v$ such that $v^t C_X v
= -1$, namely $v=(v_x)_{x \in X}$ where the integers $v_x$ are placed
at the vertices as in the following picture:
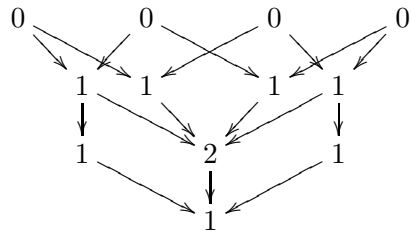
\[
\xymatrix@=1pc{0 \ar[dr] \ar[drr] & & 0 \ar[dl] \ar[drr] & & 0
\ar[dr] \ar[dll] & & 0 \ar[dl] \ar[dll]
\\ & 1 \ar[d] \ar[drr] & 1 \ar[dr] & & 1 \ar[dl] & 1 \ar[d] \ar[dll]\\
 & 1 \ar[drr] & & 2 \ar[d] & & 1 \ar[dll] \\ & & & 1 \\}
\]

On the other hand, the Coxeter matrices of the quivers $A_3$ and
$D_4$ are periodic, their orders are 4 and 6 respectively. By
Corollary~\ref{c:coxprod}, the Coxeter matrix of $X$ is periodic of
order 12.

Contrary to Example~\ref{ex:periodic}, one can show that not only the
image $\Phi_X$ of the Auslander-Reiten translation $\tau_X : \dX \to
\dX$ in the Grothendieck group is periodic, but also that actually
$\tau_X^e \simeq [d]$ for some integers $d,e \geq 1$.
\end{exmp}


\providecommand{\bysame}{\leavevmode\hbox to3em{\hrulefill}\thinspace}
\providecommand{\MR}{\relax\ifhmode\unskip\space\fi MR }
\providecommand{\MRhref}[2]{%
  \href{http://www.ams.org/mathscinet-getitem?mr=#1}{#2}
} \providecommand{\href}[2]{#2}

\end{document}